\documentclass[11pt]{article}
\usepackage{amsmath,amssymb,amsfonts,amsthm}

\newcommand{\Hom}{\mbox{Hom}\,}
\newcommand{\Ext}{\mbox{Ext}\,}

\newcommand{\Spec}{\mbox{Spec}\,}

\newcommand{\Supp}{\mbox{Supp}\,}
\newcommand{\gr}{\mbox{grade}\,}

\renewcommand{\dim}{\mbox{dim}\,}

\newcommand{\f}{\mbox{f}\,}
\newcommand{\cf}{\mbox{cf}\,}

\newcommand{\Min}{\mbox{Min}\,}

\renewcommand{\H}{\mbox{H}}
\newcommand{\V}{\mbox{V}}

\newcommand{\fa}{\frak{a}}

\newcommand{\fm}{\frak{m}}

\begin{document}
\date{}
\title{\bf Finiteness of extension functors of local cohomology modules 
\footnotetext{Email
addresses: dibaeimt@ipm.ir and yassemi@ipm.ir} }

\author{Mohammad T. Dibaei $^{\it ab}$ and Siamak Yassemi$^{\it
ac}$\\
{\small\it $(a)$ Institute for Studies in Theoretical Physics and
Mathematics, Tehran, Iran}\\
{\small\it $(b)$ Department of Mathematics, Teacher Training
University, Tehran, Iran}\\
{\small\it $(c)$ Department of Mathematics, University of Tehran,
Tehran, Iran}} \maketitle

\begin{abstract}

\noindent Let $R$ be a commutative Noetherian ring, $\fa$ an ideal
of $R$ and $M$ a finitely generated $R$--module. Let $t$ be a
non-negative integer such that $\H^i_\fa(M)$ is $\fa$--cofinite for
all $i<t$. It is well--known that $\Hom_R(R/\fa,\H^t_\fa(M))$ is
finitely generated $R$--module. In this paper we study the
finiteness of $\Ext^1_R(R/\fa,\H^t_\fa(M))$ and
$\Ext^2_R(R/\fa,\H^t_\fa(M))$.

\end{abstract}


\vspace{.3in}

\noindent{\bf 1. Introduction}

\vspace{.2in} Throughout this paper, $R$ is commutative Noetherian
ring and $\fa$ is an ideal of $R$. An $R$--module $M$ is called
$\fa$--cofinite if:

\begin{itemize}

\item[(i)] $\Supp(M)\subseteq\V(\fa)$

\item[(ii)] $\Ext^i_R(R/\fa,M)$ is finite (i.e. finitely generated)
$R$--module for all $i\ge 0$.

\end{itemize}

In [{\bf G}] Grothendieck conjectured that ``for any finite
$R$--module $M$, $\Hom_R(R/\fa,\H^i_\fa(M))$ is finite for all
$i$''. Although, Hartshorne disproved Grothendieck's conjecture (cf.
[{\bf H}]) but there are some partial answers to Grothendieck's
conjecture. For example, in [{\bf DY}, Theorem 2.1] we showed that
for a finite $R$--module $M$ and for a non--negative integer $t$ if
$\H^i_\fa(M)$ is $\fa$--cofinite for all $i<t$ then
$\Hom_R(R/\fa,\H^t_\fa(M))$ is finite.

Now it is natural to ask about the finiteness of
$\Ext^i_R(R/\fa,\H^t_\fa(M))$ for $i>0$. The first main result is to
give a partial answer for the case $i=1$, see Theorem A.

For the case $i=2$, Assadolahi and Schenzel used the spectral
sequence method to show that over a local ring $(R,\fm)$ if $M$ is a
Cohen--Macaulay $R$--module and $t=\gr(\fa, M)$ then
$\Ext^2_R(R/\fa, \H^t_\fa(M))$ is finite if and only if
$\Hom_R(R/\fa,\H^{t+1}_\fa(M))$ is so. The second main result of
this paper is to give a generalization of [{\bf AS}, Theorem 1.2]
without using spectral sequence, see Theorem B.

\vspace{.2in}

\noindent{\bf 2. Main results}

\vspace{.2in}

\noindent{\bf Theorem A:} (Finiteness of $\Ext^1_R(R/\fa,
\H^t_\fa(M))$). Let $t$ be a non--negative integer. Let $M$ be an
$R$--module such that $\Ext^{t+1}_R(R/\fa, M)$ is a finite
$R$--module (for example $M$ might be finite). If $\H^i_\fa(M)$ is
$\fa$--cofinite for all $i<t$, then
$\Ext^1_R(R/\fa,\H^{t}_\fa(M))$ is finite.

\vspace{.2in}

\noindent{\bf Theorem B:} (Finiteness of $\Ext^2_R(R/\fa,
\H^t_\fa(M))$). Let $M$ be an $R$--module such that
$\Ext^i_R(R/\fa, M)$ is finite for all $i\ge 0$ (for example $M$
might be finite). Let $t$ be a non--negative integer such that
$\H^i_\fa(M)$ is $\fa$--cofinite for all $i<t$. Then the following
statements are equivalent.

\begin{itemize}

\item[(a)] $\Hom_R(R/\fa, \H^{t+1}_\fa(M))$ is finite.

\item[(b)] If $\Ext^2_R(R/\fa, \H^t_\fa(M))$ is finite.

\end{itemize}

\vspace{.2in}

We first bring the following remark which is crucial in our proofs.

\vspace{.1in}

\noindent{\bf Remark 2.1.} Let $M$ be an $R$--module and let $E$ be
the injective hull of the $R$--module $M/\Gamma_\fa(M)$. Let
$N=E/(M/\Gamma_\fa(M))$. Then it is easy to see that the modules
$\Gamma_\fa(E)$ and $\Hom_R(R/\fa, E)$ are zero. Also from the exact
sequence $$0\to M/\Gamma_\fa(M)\to E\to N\to 0,$$ we have
$\H^i_\fa(N)\cong\H^{i+1}_\fa(M)$ and $\Ext^i_R(R/\fa,
N)\cong\Ext^{i+1}_R(R/\fa, M/\Gamma_\fa(M))$ for all $i\ge 0$. In
addition, note that $\Hom_R(R/\fa, \Gamma_\fa(N))=\Hom_R(R/\fa, N)$.

\vspace{.1in}

\noindent{\it Proof of Theorem A.} We use induction on $t$. Let
$t=0$. The short exact sequence
\begin{equation}
0\to\Gamma_\fa(M)\to M\to M/\Gamma_\fa(M)\to 0 \label{1}
\end{equation}
induces the following exact sequence
$$0=\Hom_R(R/\fa,
M/\Gamma_\fa(M))\to\Ext^1_R(R/\fa,\Gamma_\fa(M))\to\Ext^1_R(R/\fa,M),$$
and hence $\Ext^1_R(R/\fa,\Gamma_\fa(M))$ is finite.

Suppose that $t>0$ and that the case $t-1$ is settled. Since
$\Gamma_\fa(M)$ is $\fa$--cofinite, the $R$--module $\Ext^i_R(R/\fa,
\Gamma_\fa(M))$ is finite for all $i$. By the exact sequence (1),
$\Ext^{t+1}_R(R/\fa, M/\Gamma_\fa(M))$ is finite. Now by Remark 2.1
the $R$--module $\Ext^{t}_R(R/\fa, N)$ is finite and $\H^i_\fa(N)$
is $\fa$--cofinite for all $i<t-1$. Thus by induction hypothesis,
$\Ext^1_R(R/\fa, \H^{t-1}_\fa(N))$ is finite and so
$\Ext^1_R(R/\fa,\H^t_\fa(M))$ is finite.\hfill$\square$

\vspace{.2in}

\noindent{\it Proof of Theorem B.} (a)$\Rightarrow$(b) We use
induction on $t$. Let $t=0$. The short exact sequence (1) induces
the following exact sequence
$$\Ext^1_R(R/\fa,
M/\Gamma_\fa(M))\to\Ext^2_R(R/\fa,\Gamma_\fa(M))\to\Ext^2_R(R/\fa,M).$$
To show $\Ext^2_R(R/\fa, \Gamma_\fa(M))$ is finite, it is enough to
show that $\Ext^1_R(R/\fa, M/\Gamma_\fa(M))$ is finite. By Remark
2.1, we have

$$\begin{array}{ll}
\Ext^1_R(R/\fa, M/\Gamma_\fa(M)) &= \Hom_R(R/\fa, N)\\
& = \Hom_R(R/\fa, \Gamma_\fa(N)) \\
& =\Hom_R(R/\fa,\H^1_\fa(M)).
\end{array}$$
Now the assertion holds.

Suppose that $t>0$ and that the case $t-1$ is settled. Since
$\Gamma_\fa(M)$ is $\fa$--cofinite, the $R$--module $\Ext^i_R(R/\fa,
\Gamma_\fa(M))$ is finite for all $i$. Using the exact sequence (1)
we get that $\Ext^i_R(R/\fa, M/\Gamma_\fa(M))$ is finite for all
$i$. By Remark 2.1, $\Ext^i_R(R/\fa,N)$ is finite for all $i$ and
also $\Hom_R(R/\fa,\H^t_\fa(N))\cong\Hom_R(R/\fa, \H^{t+1}_\fa(M))$
is finite. By induction hypothesis the $R$--module $\Ext^2_R(R/\fa,
\H^{t-1}_\fa(N))$ is finite and hence $\Ext^2_R(R/\fa,\H^t_\fa(M))$
is finite too.

(b)$\Rightarrow$(a) We use induction on $t$. Let $t=0$. The short
exact sequence (1) induces the following exact sequence
$$\Ext^1_R(R/\fa,
M)\to\Ext^1_R(R/\fa,M/\Gamma_\fa(M))\to\Ext^2_R(R/\fa,\Gamma_\fa(M)).$$
Thus $\Ext^1_R(R/\fa, M/\Gamma_\fa(M))$ is finite. By Remark 2.1,
$\Hom_R(R/\fa, N)$ is finite and hence the $R$--module
$\Hom_R(R/\fa, \Gamma_\fa(N))$ is finite. Thus
$\Hom_R(R/\fa,\H^1_\fa(M))$ is finite.

Now let $t>0$ and that the case $t-1$ is settled. Remark 2.1 implies
that the modules $\Hom_R(R/\fa,\H^t_\fa(N))$ and $\Ext^i_R(R/\fa,
N)$ are finite for all $i$. By induction hypothesis the $R$--module
$\Ext^2_R(R/\fa,\H^{t-1}_\fa(N))$ is finite and hence
$\Ext^2_R(R/\fa,\H^t_\fa(M))$ is finite.\hfill$\square$

\vspace{.2in}

\noindent{\bf Remark 2.2.} Note that in Theorem B, one may replace
the condition ``$\Ext^i_R(R/\fa, M)$ is finite for all $i\ge 0$''
with the condition ``$\Ext^i_R(R/\fa,M)$ is finite for
$i=t+1,t+2$''.

\vspace{.2in}

Using the notation of [{\bf BS} 9.1.3], the $\fa$--finiteness
dimension of $M$ is defined as
$$ \f_\fa(M):=\Min\{j\in\Bbb N_0|\H^j_\fa(M)\,\,\,  \mbox{not
finite}\}.$$ Therefore, it is natural to define the
$\fa$--cofiniteness dimension of $M$ as
$$\cf_\fa(M):=\Min\{j\in\Bbb N_0|\H^j_\fa(M)\,\,\,  \mbox{not
$\fa$--cofinite}\}.$$ As the conventions $\cf_\fa(M)=\infty$ when
for all $j\ge
0$ the module $\H^j_\fa(M)$ is $\fa$--cofinite.\\
Using this notation we get the following corollary.

\vspace{.1in}

\noindent{\bf Corollary 2.3.} If $t\le\cf_\fa(M)$, then the
following hold:

\begin{itemize}

\item[(a)] $\Ext^i_R(R/\fa, \H^t_\fa(M))$ is
finite for $i\le 1$.

\item[(b)] $\Ext^2_R(R/\fa, \H^t_\fa(M))$ is finite if and only if
$\Hom_R(R/\fa, \H^{t+1}_\fa(M))$ is finite.

\end{itemize}

\vspace{.1in}

\noindent{\it Proof.} Part (a) follows from [{\bf DY}, Theorem
2.1] and Theorem A.

Part (b) follows from Theorem B.

\vspace{.2in}

\noindent{\bf Corollary 2.4.} Let $M$ be a finite $R$--module and
$t=\gr(\fa,M)$ or $t=\f_\fa(M)$. Then the following hold.

\begin{itemize}

\item[(a)] $\Ext^i_R(R/\fa, \H^t_\fa(M))$ is
finite for $i\le 1$.

\item[(b)] $\Ext^2_R(R/\fa,\H^t_\fa(M))$ is finite if and only if
$\Hom_R(R/\fa,\H^{t+1}_\fa(M))$ is finite.

\end{itemize}

\noindent{\it Proof.} Follows from Corollary 2.3 and the fact that
$\gr(\fa, M)\le\f_\fa(M)\le\cf_\fa(M)$.

\vspace{.2in}

The following result shows that $\Ext^2(R/\fa, \H^t_\fa(M))$ is
not always finite.

\vspace{.1in}

\noindent{\bf Corollary 2.5.} Let $(R,\fm)$ be a $3$--dimensional
analytically normal Cohen--Macaulay local domain and $\fa$ an ideal
such that $\dim R/\fa\geq 2$. If $\Spec R/\fa -\{\fm/\fa\}$ is
disconnected then $\Ext^2(R/\fa, \H^1_\fa(R))$ is not finite.

\vspace{.1in}

\noindent{\it Proof.} By [{\bf MV}, Theorem 3.9] we know that
$\Hom_R(R/\fa, \H^2_\fa(R))$ is not finite. On the other hand
$\H^0_\fa(R)$ is $\fa$--cofinite. Now the assertion follows from
Theorem B.

\vspace{.2in}

\noindent{\bf 3. Examples}

\vspace{.2in}

\noindent{\bf Example 3.1} (cf. [{\bf H}]) Let $k$ be a field and
$R=k[x,y,z,u]/(xy-zu)$. Set $\fa=(x,u)$. Then
$\Hom_R(R/\fa,\H^2_\fa(R))$ is not finite and so, by Theorem B,
$\Ext^2_R(R/\fa,\H^1_\fa(R))$ is not finite. Thus $\H^i_\fa(R)$ is
not $\fa$-cofinite for $i=1,2$.

\vspace{.2in}

\noindent{\bf Example 3.2} (cf. [{\bf MV}, Ex 3.7]) Let $k$ be a
field and $R=k[x,y,z]_{(x,y,z)}$. Set $\fa=((x)\cap(y,z))$. Then
$\Hom_R(R/\fa,\H^2_\fa(R))$ is not finite and so
$\Ext^2_R(R/\fa,\H^1_\fa(R))$ is not finite. Thus $\H^i_\fa(R)$ is
not $\fa$-cofinite for $i=1,2$.

\vspace{.2in}

\noindent{\bf Example 3.3} (cf. [{\bf MV}, Ex 3.10]) Let $k$ be a
field and $R=(k[x,y,u,v]/(xu-yv))_{(x,y,u,v)}$. Set
$\fa=(x,y)R\cap(u,v)R$. Then $\Hom_R(R/\fa,\H^2_\fa(R))$ is not
finite and so $\Ext^2_R(R/\fa,\H^1_\fa(R))$ is not finite. Thus
$\H^i_\fa(R)$ is not $\fa$-cofinite for $i=1,2$.

\vspace{.2in}

\noindent{\bf Example 3.4} (cf. [{\bf AS}, Example 4.1]) Let $k$ be
an arbitrary field. Let $R=k[|x,y,z|]$ denote the formal power
series ring in three variables. Let $\fa=(x,y)R\cap zR$. Then
$\Ext^2_R(R/\fa,\H^1_\fa(R))$ is not finite.

\vspace{.2in}

\noindent{\bf Example 3.5} (cf. [{\bf HK}, Ex 2.4]) Let $k$ be a
field of characteristic zero, and let $R=k[X_{ij}]$, for $1\le i\le
2$, $1\le j\le 3$. Let $\fa$ be the height two prime ideal of $R$
which is generated by the $2\times 2$ minors of the matrix
$(X_{ij}$. Then $\Hom_R(R/\fa,\H^3_\fa(R))$ is not finite and so
$\Ext^2_R(R/\fa,\H^2_\fa(R))$ is not finite. Thus $\H^i_\fa(R)$ is
not $\fa$-cofinite for $i=2,3$.

\vspace{.3in}

\noindent {\large\bf Acknowledgment.} This paper were initiated when
the authors were attending the 16th Algebra Seminar in Institute for
advance studies in basic sciences, Zanjan, Iran, 17--19 November
2004. It is a pleasure to thank the organizer of this seminar for
invitation and hospitality.

\baselineskip=16pt

\begin{center}
\large {\bf References}
\end{center}
\vspace{.2in}

\begin{itemize}

\item[[AS]] J. Asadollahi, P. Schenzel, {\em Some results on
associated primes of local cohomology modules}, Japan. J. Math.
(N.S.) {\bf 29} (2003), 285--296.

\item[[BS]] M. P. Brodmann, R. Y. Sharp, {\it Local Cohomology: An
Algebraic Introduction with Geometric Applications}, Cambridge
Studies in Advanced Mathematics, 60. Cambridge University Press,
Cambridge, 1998.

\item[[DY]] M. T. Dibaei, S. Yassemi, {\em Associated primes and cofiniteness
of local cohomology modules}, To appear in Manuscripta
Mathematica.

\item[[G]] A. Grothendieck, {\em Cohomologie locale des faisceaux
coh\'{e}rents et th\'{e}or\'{e}mes de Lefschetz locaux et globaux}
$(SGA$ $2)$ Advanced Studies in Pure Mathematics, Vol. 2.
North-Holland Publishing Co., Amsterdam; Masson \& Cie, Editeur,
Paris, 1968.

\item[[H]] R. Hartshorne, {\em Affine duality and cofiniteness},
Invent. Math. {\bf 9} (1970) 145--164.

\item[[HK]] C.Huneke, J. Koh, {\em Cofiniteness and vanishing of
local cohomology modules}, Math. Proc. Camb. Phil. Soc. {\bf 110}
(1991), 421--429.

\item[[MV]] T. Marley, C. Vassilev {\em Cofiniteness and associated primes of local cohomology modules},
J. Algebra {\bf 256} (2002), 180--193.

\end{itemize}

\end{document}